\documentclass[12pt]{article}
\usepackage{setspace}
\usepackage{textcomp}
\usepackage{amsmath,amsthm,amssymb}
\usepackage[left=1.5in,top=1in,right=1in,nohead]{geometry}
\usepackage{hyperref}
\usepackage{latexsym}
\usepackage{mathrsfs}
\usepackage{rawfonts}
\usepackage[dvips]{graphicx}
\usepackage[dvips]{graphicx-psmin}
\usepackage{xcolor}
\usepackage{stmaryrd}

\usepackage[sc,osf]{mathpazo}

\def\bct{\begin{center}}
\def\ect{\end{center}}
\def\beg{\begin}
\def\bit{\begin{itemize}}
\def\eit{\end{itemize}}

\def\<{\langle}
\def\>{\rangle}

\def\mbb{\mathbb}

\def\mco{\mathcal O}

\def\ni{\noindent}

\def\tn{\textnormal}
\def\wt{\widetilde}

\def\ZZ{{\mathbb Z}}
\def\CP{{\mathbb C \mathbb P}}
\DeclareMathOperator{\Fix}{Fix}

\newtheorem{thm}{Theorem}[section]
\newtheorem{lem}[thm]{Lemma}

\newtheorem{cor}[thm]{Corollary}
\newtheorem{rmk}[thm]{Remark}

\title{On the Curvature of Einstein-Hermitian Surfaces 
} \author{Mustafa Kalafat \and Caner Koca}

\begin{document}
\maketitle
\begin{abstract} We give a mathematical exposition of the Page metric, and introduce an efficient coordinate system for it. 
We carefully examine the submanifolds of the underlying 
smooth manifold, and show that the Page metric does not have positive holomorphic bisectional curvature. We exhibit a holomorphic subsurface with flat normal bundle. 
We also give another proof of the fact 
that a compact complex surface together with an Einstein-Hermitian metric of positive orthogonal bisectional curvature is biholomorphically isometric to 
the complex projective plane with its Fubini-Study metric up to rescaling. This result relaxes the K\"ahler condition in Berger's theorem, and the positivity condition on sectional curvature in a theorem proved by the second author. 

  \end{abstract}


\section{Introduction}

Let $(M,J)$ be a complex manifold. A Riemannian metric $g$ on $M$ is called {\em Hermitian} 
if the complex structure $J:TM\rightarrow TM$ is an orthogonal transformation at every point on $M$ with respect to the metric $g$, that is, 
 $g(X,Y)=g(JX,JY)$ for tangent vectors $X,Y\in T_p M$ for all $p\in M$. In this case, the triple $(M,g,J)$ is called a {\em Hermitian manifold}.
For Hermitian metrics we have further notions of curvature related to complex structure: The {\em holomorphic bisectional curvature} in the direction of a pair of unit tangent vectors $X,Y\in T_pM$ is defined as
$$\tn{H}(X,Y):=\tn{Rm}(X,JX,Y,JY).$$
If one applies the algebraic Bianchi identity, and $J$-invariance in the K\"ahler case, it is easy to see that this is the sum of sectional curvatures of the planes spanned by $X,Y$ and $X,JY$. 
We call this as the {\em summation identity} in the K\"ahler case which gives some visual insight. Bisectional curvature 
is actually \emph{not} an invariant of the plane spanned by the vectors $X,Y$. Rather, it is an 
invariant of the holomorphic planes spanned by them. As a special case, if one takes the two vectors 
identical, then the result coincides with the sectional curvature of the holomorphic plane spanned. 
This is 
called the {\em holomorphic (sectional) curvature} in that direction vector. Although  
it can be considered as a map on the sphere $S^{2n-1}(T_p M) \rightarrow \mbb R$ at each point, 
just in the case of bisectional curvature, this map is an invariant of the holomorphic plane, and therefore can be considered as a map on the complex Grassmannian of one lower real dimension since it is constant on the Hopf circle fibers. 
Positivity of the sectional curvature implies that of bisectional curvature, which 
implies positivity of holomorphic curvature. However, the converses are not necessarily 
true in general. 

In this paper, we work on some explicit 4-manifolds to understand various notions of curvature. We also  
prove a uniformization theorem for positive bisectional curvature. 
First, let us review some well-known theorems in special cases. 

\begin{thm}[Frankel conjecture, Siu-Yau Thm \cite{SY:1980frankel}] \label{bisecfrankeln}
Every compact K\" ahler manifold of positive bisectional holomorphic curvature is biholomorphic to 
the complex projective space.\end{thm}
\ni This theorem does not, however, specify the metric in question. Nevertheless, if we in addition 
assume that the metric is \emph{Einstein}, then the metric is unique, too:
\begin{thm}[\cite{berger65,goldbergkobayashi67}]  \label{bgk}
An $n$-dimensional compact connected K\" ahler manifold with an Einstein (or constant scalar 
curvature) metric of positive holomorphic bisectional curvature is globally isometric to $\mbb{CP}
_n$ with the Fubini-Study metric (up to rescaling).
\end{thm}

\ni Our aim here is to relax the K\"ahler condition on the metric in Theorem \ref{bgk} to 
merely being Hermitian in dimension 4. In this case the summation identity 
is no longer valid. Our main theorem is:

\begin{thm} \label{mainthm} If $(M,g,J)$
is a compact complex surface together with an Einstein-Hermitian metric of positive holomorphic 
bisectional curvature, then it is biholomorphically isometric to $(\mbb{CP}_2,g_{FS})$, the complex 
projective plane with its Fubini-Study metric up to rescaling. 
\end{thm}

\ni We note that an analogous result with the positivity assumption on the 
\emph{sectional} curvature is proved in \cite{kocapage} by the second author. 
However, in the Hermitian case, positivity of 
sectional and bisectional curvatures are not related in general due to invalidity of 
the summation identity. 
Consequently the techniques in that paper cannot be directly applied here.

Einstein-Hermitian metrics on compact complex surfaces are classified by LeBrun:\newpage
\begin{thm}[\cite{LeB:2011}]\label{lebruntheorem}
If $(M,g,J)$ is a compact complex surface together with an Einstein-Hermitian metric, then only one of the following holds.
\begin{enumerate}
\item $g$ is K\" ahler-Einstein(KE).

\item $M$ is biholomorphic to $\mbb{CP}_2\sharp\,\overline{\mbb{CP}}_2$ and
$g$ is the Page metric (up to rescaling).

\item $M$ is biholomorphic to $\mbb{CP}_2\sharp\, 2\overline{\mbb{CP}}_2$ and
$g$ is the Chen-LeBrun-Weber metric (up to rescaling).
\end{enumerate}
\end{thm}

\ni In other words, an Einstein-Hermitian metric is either K\"ahler-Einstein to start with, or is one of the two exceptional models. These exceptional Einstein metrics are non-K\"ahler, but they are 
\emph{conformally K\"ahler}. 
By Theorem \ref{bgk} we only need to consider the non-K\"ahler case. 
We start with the Page metric case. We give an introduction to this metric, and show that it does not 
have positive curvature everywhere. 
One of the key facts in the proof of this result is Frankel's Theorem \cite{Fra:1961} which states that totally geodesic submanifolds of complementary dimensions
on positively curved manifolds necessarily intersect. Since the
Page metric has an explicit form, we are also able to give a computational
proof of the failure of positivity. Secondly we introduce Euler coordinates, and use 
Dragomir-Grimaldi's theorem \cite{dragomirgrimaldi91} to get the analogous result on bisectional curvature. These coordinates are especially useful for Page metric and its submanifolds. We hope that they will be useful for others who are interested in local computations in 4-dimensional geometry. For another application of these coordinates see \cite{kalafatsari}.
 
On the other hand, the Chen-LeBrun-Weber metric does not have such an explicit formula. Therefore, this explicit analysis is not an available option at this time. We need to prove a more general result to handle this case. 
Using curvature estimates and Weitzenb\" ock formula techniques we prove the 
following result.

\vspace{.05in}

\noindent {\bf Theorem \ref{finalthm}.}
{\em Let $M$ be a compact Einstein-Hermitian 4-manifold of positive holomorphic bisectional curvature. 
Then the Betti number $b^2_-$ vanishes.}

\vspace{.05in}  

\ni As a consequence we have, 

\begin{cor}
There is no Einstein-Hermitian metric of positive holomorphic bisectional curvature on any blow up 
of $\mbb{CP}_2$.
\end{cor}

\ni Since the underlying smooth 4-manifolds of the exceptional cases $\mbb{CP}_2\sharp
\,\overline{\mbb{CP}}_2$ and $\mbb{CP}_2\sharp\, 2\overline{\mbb{CP}}_2$ have non-zero $b^2_-$, we 
deduce the following corollary. 
\begin{cor}
The Page and Chen-LeBrun-Weber metrics are not of positive holomorphic bisectional curvature.
\end{cor}

\ni This eliminates the later two possibilities in LeBrun's Classification Theorem \ref{lebruntheorem}. In the remaining K\" ahler-Einstein case, we apply  the Berger -- Goldberg-Kobayashi Theorem \ref{bgk}, and the proof of our Theorem \ref{mainthm} follows. 

We note that the authors proved a more general result in a different article \cite{confk} in the conformally K\" ahler case, which uses slightly different Weitzenb\" ock techniques. Besides that, the discussion on Page metric and Euler coordinate computations are the main and completely new material here. 
We give a better understanding of this important metric to the reader in this paper.

 In \S\ref{pagemath} we give a careful topological analysis of the Page metric and show that it 
 does  not have positive sectional curvature everywhere. In \S\ref{secpagesubmanifolds} we provide Euler coordinates and give an alternative proof of the fact that it the Page metric is not of positive 
 holomorphic bisectional curvature by exhibiting a subsurface with flat normal bundle. In \S\ref{secestimates} we prove some estimates and 
classify Einstein-Hermitian 4-manifolds of positive bisectional curvature. 
 
\noindent{\bf Acknowledgements.} We thank Claude LeBrun for his encouragement. 
Also thanks to \"O. Kelek\c ci and C. Tezer  
 for useful discussions. 
 This work is partially supported 
by the grant $\sharp$113F159 of T\"ubitak (Turkish science and research council).



\section{Page metric}\label{pagemath}
In this section we give a rigorous mathematical exposition of the Page metric, and describe its topology in detail. At the end we prove that it does not have positive sectional curvature everywhere. This section and its figures are part of C. Koca's thesis \cite{kocathesis}. 

The Page metric was discovered by D. Page in 1978 as a limiting metric of Kerr-de Sitter solution (see \cite{P:1978}).  It is the unique Einstein-Hermitian non K\"ahler metric on the blow up of complex projective plane. To define it formally, we first think of the following metric on the product $S^3\times I$ where I is the closed interval $[0,\pi]$:
\[
g=V(r) dr^2 + f(r) (\sigma_1^2 +\sigma_2^2) + \frac{C \sin^2 r}{V(r)} \sigma_3^2
\]
where the coefficient functions are given by the following expressions 
\footnote{In the original paper \cite{P:1978} of Page, there is a typo in the equation (41). The last term inside the curly brackets has to be divided by 4 
and $(3-\nu^2)^2$ replaced by $(3+\nu^2)^2$. 
See \cite{page2arxiv} for the information.} 
$$V(r)  = \frac{1-a^2 \cos^2 r}{3-a^2-a^2 (1+a^2) \cos^2 r}$$
$$f(r)  = 4\,\frac{1-a^2 cos^2 r}{3+6a^2-a^4}$$
$$C  = \frac{1}{(3+a^2)^2}$$

\begin{figure}[!h]
\begin{center}~\\ \vspace{16mm}
\includegraphics[scale=1.2]{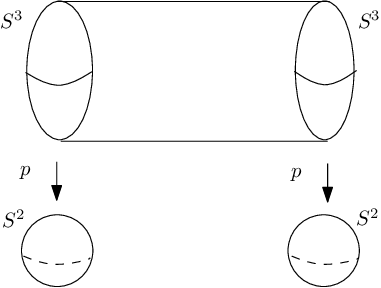}
\end{center}
\caption{\em The manifold $S^3\times I$ with its two boundary components.}\label{figur1}
\end{figure}

\begin{figure}
\begin{center}
\includegraphics[scale=1]{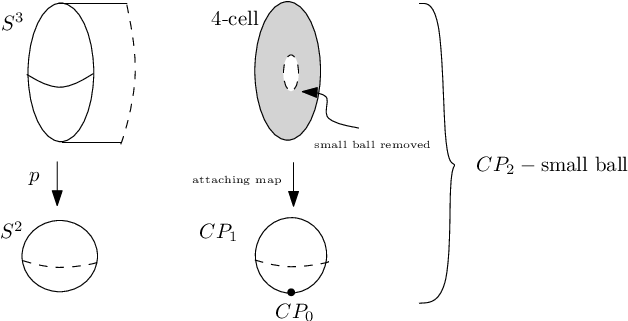}
\end{center}
\caption{\em Boundary identifications to produce the closed 4-manifold $\CP_2\sharp\,\overline{\CP}_2$.}
\label{figur2}
\end{figure}

\ni and $a$ is the unique positive root of $a^4+4a^3-6a^2+12a-3=0$. Here, $\sigma_1,\sigma_2,\sigma_3$ is the standard left invariant 1-forms on the Lie group $SU(2)\approx S^3$. 
At the endpoints $r=0$ and $\pi$, we see from the formula that the metric shrinks to a round metric on $S^2$. 
Thus, $g$ descends to a metric, denoted by $g_{\mathrm{Page}}$, on the quotient $(S^3\times I) 
\diagup \sim$ where $\sim$ identifies the fibers of the Hopf fibration $p:S^3 \rightarrow S^2$ on 
the two ends $S^3\times \{0\}$ and $S^3\times\{\pi\}$ of the cylinder $S^3\times I$. See Figure 
\ref{figur1}. 
\noindent 
The resulting manifold is indeed the connected sum $\CP_2\sharp\,\overline{\CP}_2$. To see this, recall that in the cell decomposition of $\CP_2$, the attaching map from the boundary of the 4-cell (which is $S^3$) to the 2-skeleton (which is $\CP_1\approx S^2$) is given by the Hopf map \cite{Hat:2002}. So, if we cut the cylinder $S^3\times I$ in two halves and identify the Hopf fibers of $S^3$ at each end, we get $\CP_2 - \textnormal{\{small ball\}}$. Since the right and left halves have different orientations, we obtain $\CP_2\sharp\,\overline{\CP}_2$ in the quotient. See Figure 
\ref{figur2} for assistance. 
\noindent Next, we will prove that the Page metric is not of positive sectional curvature. We will use the following classical theorem by Frankel.
\begin{thm}[\cite{Fra:1961}]
Let $M$ be a smooth $n$-manifold, and let $g$ be a complete Riemannian metric of positive sectional curvature. If $X$ and $Y$ are two compact totally geodesic submanifolds of dimensions $d_1$ and $d_2$ such that $d_1+d_2\geq n$, then $X$ and $Y$ intersect.
\end{thm}
In our case, the two 2-spheres on each end of the above quotient will play the role of $X$ and $Y$. They are compact and the dimensions add up to $4$. So it remains to show that those two submanifolds are totally geodesic with respect to $g_{\mathrm{Page}}$. Since they are obviously disjoint, this will imply that $g_{\mathrm{Page}}$ cannot have positive sectional curvature. 
There is a very well-known lemma to detect totally geodesic submanifolds:

\begin{lem}\label{lemma}
Let $(M,g)$ be a Riemannian manifold. If $f$ is an isometry, then each connected component of the fix point set $\Fix(f)$ of $f$ is a totally geodesic submanifold of $M$.
\end{lem}

So, below we will show that there is an isometry of the Page metric whose fix point set is precisely the two end spheres. What are the isometries of the Page metric? Derdzi\' nski showed that the Page metric is 
indeed conformal to one of Calabi's extremal K\"ahler metrics on $\CP_2\sharp\,\overline{\CP}_2$ in \cite{Derdzinski}. On 
the other hand, the identity component of the isometry group of extremal K\"ahler metrics is a 
maximal compact subgroup of the identity component of the automorphism group \cite{C2:1985}. In the case 
of $\CP_2\sharp\,\overline{\CP}_2$, this implies that the identity component of the isometry group of 
the Page metric is $U(2) = (SU(2)\times S^1)/\ZZ_2$. By the formula of the metric, we see that the 
isometries in the $SU(2)$ component are precisely given by the left multiplication action of $SU(2)
$ on the first factor of $S^3\times I$. Note that the forms $\sigma_i$, $i=1,2,3$ are invariant 
under the action, but the action on the 3-spheres $S^3\times\{r\}$, $r\in(0,\pi)$ is fixed-point-
free! The metric is invariant under this action as the coefficients of the metric only depend on 
the parameter $r$.

Now, let us see what happens at the endpoints $r=0$ and $r=\pi$: It is well-known that the action of $U\in SU(2)$ on the 2-sphere $S^2$ (after the quotient) is given by the conjugation $A\mapsto U AU^{-1}$, where we regard the $2\times 2$ complex matrix  $A=x\sigma_1+y\sigma_2 +z\sigma_3$ with $x^2+y^2+z^2=1$ as a point of $S^2$. It is now straightforward to see that the action of $-I\in SU(2)$ is trivial on $S^2$ (since $(-I)A(-I)^{-1}=A)$; thus, it fixes every point on $S^2$. Therefore, we conclude that the fixed point set of the isometry given by the ``antipodal map" $-I\in SU(2)$ consists of the two $2$-spheres at each end of the quotient $\left((S^3\times I)\diagup\sim\right)\approx\CP_2\sharp\,\overline{\CP}_2$. Note that, indeed, there is an $S^1$-family of isometries generated by rotation in direction of $\sigma_3$ having the exact same fixed point set. 
So we showed that there are two disjoint compact totally geodesic submanifolds of $\CP_2\sharp\,
\overline{\CP}_2$. Therefore, 
Frankel's theorem implies the following. 

\beg{thm} The sectional curvature of the Page metric is not everywhere positive. 
\end{thm}

Finally, we note that we can actually show the failure of positivity \emph{directly by brute-force} using tensor calculus: Introduce a new coordinate function $x:=\cos(r)$, so that the metric becomes
\[
g=U^2(x) dx^2 + g^2(x) (\sigma_1^2 +\sigma_2^2) + \frac{D^2}{W(x)} \sigma_3^2
\]
where the coefficient functions are given as

$$U(x)  = \sqrt{\frac{1-a^2 x^2}{(3-a^2-a^2 (1+a^2) x^2)(1-x^2)}}$$

$$g(x)  = 2\,\sqrt{\frac{1-a^2 x^2}{3+6a^2-a^4}}$$

$$D   = \frac{1}{3+a^2}$$

\ni and choose the following vierbein: $\{e^0,e^1,e^2,e^3\} := \{U dx, g\sigma_1, g\sigma_2, DU^{-1}\sigma_3\}$. Then by a standard tensor calculus, we see that the sectional curvature of the plane generated by $e_0$ and $e_1$ is given by $$K_{01}=2\,\frac{g'U'-g'' U}{gU^3}.$$ Using a computer program like \emph{Maple}, one can easily verify that this function $K_{01}(x)$ can take both positive and negative values for $x\in(-1,1)$.

\section{Euler coordinates 
and Flat bundles}\label{secpagesubmanifolds}

In this section we will introduce an efficient coordinate system and use it to show that the Page 
metric is not of positive holomorphic bisectional curvature by explicitly analyzing the  
submanifolds of the underlying smooth manifold. 
We would like to use the following theorem of Dragomir and Grimaldi. 
See the book \cite{dragomirornea} on locally conformal K\" ahler (l.c.K.) geometry p.157 for an exposition.
\beg{thm}[\cite{dragomirgrimaldi91}] \label{dragomirgrimaldi91}
Let $S$ be a complex submanifold of the l.c.K. manifold $M$.
If $M$ has positive holomorphic bisectional curvature everywhere,
then the normal bundle of the given immersion $S\subset M$ admits no parallel sections.
\end{thm}
In order to make use of this theorem, we need to analyze the complex submanifolds of $\mbb{CP}_2\sharp\,\overline{\mbb{CP}}_2$. For this purpose we
use {\em Euler angles} \cite{zhangeulerforms} 
 on the $S^3\subset \mbb R^4$ which e.g. realizes the Hopf fibration in the best.
$$0\leq\theta \leq \pi, ~~0\leq\phi\leq 2\pi, ~~0\leq{{\psi}}\leq 4\pi $$
\vspace{-.5cm}
$$x_1:=r\cos{\theta\over 2}\cos{\psi+\phi\over 2}$$
$$x_2:=~~~~''~~~~ 
\sin{\psi+\phi\over 2}$$
$$x_3:=r\sin{\theta\over 2}\cos{\psi-\phi\over 2}$$
$$x_0:=~~~~''~~~~ 
\sin{\psi-\phi\over 2}$$
where the {\em Hopf fibration} in these coordinates is just a projection \cite{fip},
$$h:S^3\longrightarrow S^2, ~~~~~~~~h(\theta,{{\psi}},\phi)=(-\phi,\theta).$$
Here the exchange $\phi\longleftrightarrow\theta$ is needed to relate to the calculus angles on $S^2$. Changing ${{\psi}}$ does not change the element in the image.  So, whenever the image $\phi,\theta$ is fixed, ${{\psi}}$ parametrizes the Hopf circle (fiber).

\vspace{3mm}

\ni An orthonormal, invariant coframe $\{\sigma_1,\sigma_2,\sigma_3\}$ on $S^3$ is given as follows:   
$$\beg{array}{l}
\sigma_1= (x_1dx_0-x_0dx_1+x_2dx_3-x_3dx_2)/r^2=
( \sin\psi d\theta-\sin\theta\cos\psi d\phi)/2 \nonumber\\ [4\jot] 
\sigma_2= (x_2dx_0-x_0dx_2+x_3dx_1-x_1dx_3)/r^2=
(-\cos\psi d\theta-\sin\theta\sin\psi d\phi)/2 \nonumber \\ [4\jot] 
\sigma_3= (x_3dx_0-x_0dx_3+x_1dx_2-x_2dx_1)/r^2=
(d\psi+\cos\theta d\phi)/2. \nonumber \label{relation5} 
\end{array}$$

\ni One can check in a straightforward manner the identities,
$$d\sigma_1=2\sigma_2\wedge\sigma_3 ~~~~\tn{and}~~~~ \sigma_1^2+\sigma_2^2=(d\theta^2+\sin^2\theta d\phi^2)/4.$$
Plugging these into the Page metric's expression we get
$$\hspace{-20mm}g_{\mathrm{Page}}=Vdr^2+\left\{ {f\over 4}\sin^2\theta +{C\sin^2r \cos^2\theta\over 4V(r)} \right\}\hspace{-1mm}d\phi^2+ {C\sin^2r\over 4V(r)}d\psi^2
+ {C\sin^2r \cos\theta\over 4V(r)}
(d\psi\otimes d\phi+d\phi\otimes d\psi)+{f\over 4}d\theta^2.$$
Letting~ $U:=\sqrt{V(r)}, ~h:=\sqrt{f}, ~D:=\sqrt C$ we have the  Vierbein i.e. orthonormal coframe 


\begin{multline*}\{e^0,e^1,e^2,e^3\}=\left\{ ~Udr,
~  {Dh\over 2\sqrt{C \cot^2\theta+Vf\csc^2r}}\, d\psi,\right. \\ \left. 
{D\sin r\over 2U\sqrt{1+C^{-1}Vf\tan^2\theta\sin^{-2}r} }\, 
(d\psi+(C^{-1}Vf\tan\theta\sin\theta\sin^{-2}r+\cos\theta)\,d\phi),
~{h\over 2}\,d\theta~ \right\}  \end{multline*}
 
Now, we are in a position to analyze some of the subsurfaces easily \cite{geofrust,anderssonthesis}. For example keeping $r_0,\theta_0$ fixed and varying $\psi,\phi$, one obtains tori. See more on the subsurfaces of the Page space at \cite{kalafatsari}. 
We are interested in complex submanifolds. For this purpose,
this time keep $\phi_0,\theta_0$ fixed, vary $r,\psi$ to obtain complex spheres i.e. rational curves as follows.  This captures a series of nearby Hopf fibers corresponding to $r\in[0,\pi]$ hence a cylinder
$S^1\times I$ inside $S^3\times I$ which projects to a sphere under Hopf identification of circles at the two ends. These are complex submanifolds, since they correspond to fibers coming from Hirzebruch Surface description/fibration. This follows easily by recalling that the complex line bundle $\mathcal O(1)$ over $\mbb{CP}_1$ can also be described as the quotient
$$(S^3\times \mathbb C)/\propto \quad\textrm{where}\quad (x,z)\propto(\lambda x,\lambda z),\quad \lambda\in S^1\subset\mathbb C.$$
There is an obvious diffeomorphism between
$${{(S^3\times [0,\pi))}\over \sim}\approx{{(S^3\times [0,\infty))}\over\sim} \quad \textrm{and}\quad \mathcal O(1)=(S^3\times \mathbb C)/\propto;$$ and under this map the quotient of the cylinder in the previous paragraph is mapped to the complex fiber of the complex line bundle $\mathcal O(1)$ over the point $(\phi_0,\theta_0)\in\mbb{CP}_1$. Adding the remaining point at infinity shows that the quotient of the cylinder in  $(S^3\times I)/\sim$ corresponds to the Hirzebruch fiber of the first Hirzebruch surface $\mbb
P(\mco\oplus\mco(1))\rightarrow\mbb{CP}_1$. These fibers are obvously complex rational curves. The Hopf spheres at the two ends of $S^3\times I$ correspond to $0,\infty$-sections of this fibration
$$\begin{array}{ccc@{}c}
\mbb{CP}_1 & \to & \mbb{CP}_2\sharp\,\overline{\mbb{CP}}_2\hspace{-1mm}= & \,\mbb
P(\mco\oplus\mco(1)) \\
&&& \downarrow   \\
&&& \mbb{CP}_1
\end{array}$$

Now let's find the curvature of the normal bundle of these fibers.
Finding nontrivial connection 1-forms of the Page metric with respect to this basis 
is not an easy task. We refer the reader to \cite{kalafatsari} for this type of approach with a different choice of vierbein. 
Instead we will compute the Christoffel symbols and the coefficients of the 
Riemann curvature tensor to figure out the curvature 2-form
of the normal plane bundle of the spheres that we are working on. 
Here the submanifold directions  are  $0,1$ and normal bundle directions are $2,3$.
Since $\theta=\theta_0$ is constant (together with $\phi_0$),
using one of the computer algebra systems the {curvature 2-form} of the normal bundle $\tn{N} \mbb{CP}_1$ can be computed as \cite{egh},
{\large
$${
\renewcommand{\arraystretch}{2} 
\begin{array}{rcl} 
\tilde{R}^2_{\phantom{1} 3} 
& = & d\tilde{\omega}^2_{\phantom{1} 3}+\tilde{\omega}^2_{\phantom{1} 1}\wedge\tilde{\omega}^1_{\phantom{1} 3} \\
& = & {1\over 2}\tilde{R}^2_{\phantom{1} 3 cd}\, e^c\wedge e^d\\
& = &  {D ( V f' + f ( V' - 2V \cot r )) \sin r \tan\theta \over 
      2 f^{3/2} V^{3/2} \sqrt{1 + C^{-1}Vf\csc^2 r\tan^2\theta} }\,
      {\Huge d r\wedge d\psi.}   
\end{array} }$$  }

\ni Here recall that $V$ and $f$ are functions of $r$. 
If we focus on one of the spheres where  $\theta=0$ and $\phi=\phi_0$, this curvature 2-form vanishes. 
So that the normal bundle of this type of sphere is flat. Then the parallel translation on this sphere depends only on the homotopy class, which is unique because of simple-connectivity. Hence, parallel translation is totally path independent. Now starting with two linearly independent vectors of the normal bundle at a point, one extends them by parallel translation to the whole sphere. Since lengths are preserved during the process, these extensions are nowhere vanishing. So we obtain two parallel sections of the normal bundle. Two nowhere zero sections of a plane bundle trivializes it, so we have a trivial, flat normal bundle.
We are ready to state the main result of this section.

\beg{thm} The holomorphic bisectional curvature of the Page metric on
$\mbb{CP}_2\sharp\,\overline{\mbb{CP}}_2$ is not everywhere positive. \end{thm}

\beg{proof} The spheres above are complex submanifolds and the Page metric is conformally K\" ahler so certainly l.c.K. Their normal bundle has nontrivial parallel sections. Therefore, we can apply 
Theorem \ref{dragomirgrimaldi91}. \end{proof}

\section{Estimates and the bisectional curvature 
}\label{secestimates}

We first describe the 2-form interpretation of the planes in the tangent space which will be very useful. 
Sectional curvatures at a point $p\in M$ can be thought as a function on the Grassmannian of oriented two planes in the corresponding tangent space.
 $$\tn{sec} : G_2^+(T_p M) \longrightarrow \mbb R$$

\noindent In dimension $4$ we have a nice description of this Grassmannian in terms of forms
 $$G_2^+(\mbb R^4)\approx 
 \{(\alpha,\beta)\in \Lambda^2_+\oplus\Lambda^2_- : |\alpha|=|\beta|=1/\sqrt{2}\}\approx S^2\times S^2.$$

\ni See \cite{jeffdg} as a reference. Starting with a plane $\pi$, one can choose a special orthogonal basis $\{e_1,e_2\}$ which corresponds to the form $\sigma=e^1\wedge e^2\in\Lambda^2$ using metric duals. One can choose the basis 
in a unique way if the following conditions, 
$$e^1\wedge e^2=\alpha+\beta ~~~ \tn{for} ~~~ \alpha\in\Lambda^2_+,~ \beta\in \Lambda^2_-,~|\alpha|=|\beta|=1/\sqrt{2}$$ 
are imposed. 
Conversely, starting with a decomposable 2-form $\omega\in\Lambda^2$ i.e. $\omega=\theta\wedge\delta$ for some $\theta,\delta\in \Lambda^1$. 
The duals  $\{\theta^\sharp,\delta^\sharp\}$ gives an oriented basis for a plane.  
In this correspondence, a complex plane correspond to a form in the form 
${\omega\over 2} + \varphi$  for an anti-self-dual 2-form $\varphi$. 
If a plane $\tilde{\sigma}$ 
corresponds to 
$(\alpha,\beta)$ or simply $\sigma=\alpha+\beta \in \Lambda^2$, 
$$\tn{sec}(\tilde{\sigma}
)=\tn{Rm}(\sigma,\sigma)=\langle \mathcal R(\sigma), \sigma \rangle$$ 
where $\mathcal R:\Lambda^2\to \Lambda^2$ is the {curvature operator.} 
Recall that for complex planes $\tilde{\sigma},\tilde{\tau}$  
we compute the {\em bisectional curvature} by 
$$H(\tilde{\sigma},\tilde{\tau})=\tn{Rm}(\sigma,\tau).$$


Let us recall the decomposition of the curvature operator $\mathcal R :\Lambda^2\rightarrow \Lambda^2$. If $(M,g)$ is any oriented 4-manifold, then the decomposition $\Lambda^2 =\Lambda^2_+ \oplus \Lambda^2_-$ implies that the curvature operator $\mathcal R$ can be decomposed as

$${\mathcal R}=
\left(
\mbox{
\begin{tabular}{c|c}
&\\
$W^++\frac{s}{12}$&$\mathring{r}$\\ &\\
\cline{1-2}&\\
$\mathring{r}$ & $W^-+\frac{s}{12}$\\&\\
\end{tabular}
} \right)$$  

\ni where $\Lambda^2_{\pm}$ stands for the self-dual and anti-self-dual $2$-forms, i.e. $\Lambda^2_{\pm} = \{\phi\in\Lambda^2 : *\phi = \pm \phi \}$, where $*$ is the Hodge-$*$ operator determined by the metric $g$. 
$W^{\pm}$ is the self-dual/anti-self-dual Weyl curvature tensor, $s$ is the scalar curvature and $\mathring r$ stands for the trace-free part of the Ricci curvature tensor $r$.

Now we are ready to state and prove our first estimate lemma. Let $(M,J,g)$ be a compact Einstein-Hermitian 4-manifold. In this case, according to a theorem of LeBrun in \cite{LeB:1995}, such a space is either K\" ahler, or else they are conformally related to a K\" ahler metric $\tilde g$ with positive scalar curvature $\tilde s$ in such a way that $g=\tilde s^{-2} \tilde g$. Moreover, the non-K\"ahler Einstein-Hermitian metrics of Theorem \ref{lebruntheorem}, namely the Page metric and the Chen-LeBrun-Weber metric, are known to have (constant) positive scalar curvature. From now on, we will decorate the curvature expressions of the conformally related K\"ahler metric $\tilde g$ with tilde (e.g. $\tilde s, \tilde R, \tilde W$ etc.). We will denote the inner product of tensors with respect to $\tilde g$ by $\langle , \rangle_{\tilde g}$.

\beg{lem}[First Estimate] \label{firstestimate} Let $(M,J,g)$ be an Einstein-Hermitian metric which is not K\"ahler. Let $\tilde g$ be the conformally related K\"ahler metric such that $g=\tilde s^{-2} \tilde g$. If $g$ has positive holomorphic bisectional curvature  and $\lambda$ is an eigenvalue of the operator $\tilde W_-:\Lambda^2_- \rightarrow \Lambda^2_-$, then we have
$$\lambda<{\tilde{s}\over 6}.$$
\end{lem}

\beg{proof} Take two arbitrary unit tangent vector $X,Y$ in $T_p M$ at an arbitrary point $p$, and let  $\varphi,\psi\in\Lambda^2_+$ be the anti-self-dual 2-forms of length $1/\sqrt{2}$ corresponding to the complex lines $\{X,JX\}$ and $\{Y,JY\}$, as described above. Keeping in mind that $\mathring r = 0$ since $g$ is Einstein and $\Lambda^2_+$ and $\Lambda^2_-$ are orthogonal with respect to $g$, we see by the decomposition of $\mathcal R$ that
{\large $$\beg{array}{rcl}
\tn{H}(X,Y)
&=& \langle\mathcal R\left( {\omega\over 2}+\varphi\right)\, , \, {\omega\over 2}+\psi\rangle\\[10pt]
&=&\langle\left( W_++{s\over 12}\right){\omega\over 2}+
\left( W_-+{s\over 12}\right)\varphi , {\omega\over 2}+\psi\rangle\\[10pt]
&=& {1\over 4}\langle W_+ \omega,\omega\rangle + {s\over 48}\langle \omega,\omega\rangle + \langle W_- \varphi,\psi\rangle + {s\over 12}\langle\varphi,\psi\rangle\\[10pt]
&=&{1\over 4}\langle W_+ \omega,\omega\rangle
+{s\over 12}\left\{{1\over 2}+\langle\varphi,\psi\rangle\right\}
+\langle W_-\varphi,\psi\rangle.
\end{array}$$}
Now, since $g=\tilde s^{-2}\tilde g$, we see that
\begin{align*}
\langle W_+ \omega,\omega\rangle &= {(W_+)_{ij}}^{kl}\, \omega_{kl}\, \omega^{ij} = \tilde s^2 {(\tilde W_-)_{ij}}^{kl}\,\, \tilde s^{-2}\tilde\omega_{kl}\,\, \tilde s^{2}\tilde\omega^{ij} = \tilde s^2 \langle \tilde W_+ \tilde\omega,\tilde\omega\rangle_{\tilde g}
\end{align*}
It is well known \cite{Derdzinski} that for K\"ahler metrics the K\"ahler form is an eigenvector of the self-dual Weyl operator, more explicitly $\tilde W_+ \tilde \omega= {\tilde s\over 6} \tilde \omega$. Thus we see that $\langle \tilde W_+ \tilde\omega,\tilde\omega\rangle_{\tilde g} = {\tilde s \over 6 }\langle \tilde \omega,\tilde\omega\rangle_{\tilde g} = {\tilde s\over 3} $  and hence
$$
{\tn H}(X,Y)={\tilde{s}^3\over 12}
+{s\over 12}\left\{{1\over 2}+\langle\varphi,\psi\rangle\right\}
+\langle W_-\varphi,\psi\rangle
$$
Now, since the operator $\tilde W_-:\Lambda^2_-\rightarrow\Lambda^2_-$ is symmetric, it is diagonalizable. Let $\lambda$ be an eigenvalue, $\varphi\in{\tilde E}_\lambda\subset \Lambda^2_-$ be a corresponding eigenvector with $|\varphi|=1/\sqrt{2}$ (norm taken with respect to $g$). Choosing $\psi:=-\varphi$ in the above equation yields,
{\large
$$\beg{array}{rcl}
\tn{H}(X,Y)
&=&{\tilde{s}^3\over 12}
+{s\over 12}\left\{{1\over 2}-|\varphi|^2 \right\}
-\langle \tilde s^2 \tilde W \varphi, \varphi\rangle\\[10pt]
&=& {{\tilde s^3}\over 12} + \tilde s^2\langle  -\lambda \varphi, \varphi\rangle
={{\tilde s^3}\over 12} - \lambda \tilde s^2|\varphi|^2\\[10pt]
&=& {\tilde{s}^2\over 2} \left\{  {\tilde{s}\over 6}
-\lambda \right\} >0.
\end{array}
$$}
and hence ${\tilde{s}\over 6} -\lambda >0$, as required.
\end{proof}

\beg{lem}[Second Estimate] \label{secondestimate} Let $(M,J,g)$ be an Einstein-Hermitian metric which is not K\"ahler.  If $g$ has positive holomorphic bisectional curvature, then for all $\varphi\in\Lambda_-^2$ with $|\varphi|=1/\sqrt{2}$, we have
{ $${\tilde{s}\over 12}-\langle {\tilde W}_-\varphi,\varphi \rangle >0.$$} \end{lem}

\beg{proof}
Inserting $\psi=-\varphi$ into the first equation of the previous proof yields,
{\large $$\beg{array}{rcl}
\tn{H}(\varphi,-\varphi)
&=&{\tilde{s}^3\over 12}
+{s\over 12}\left\{{1\over 2}-|\varphi|^2 \right\}
-\langle W_-\varphi,\varphi \rangle={\tilde{s}^3\over 12} - \langle {\tilde s^2} W_-\varphi,\varphi \rangle\\[10pt]
&=& \tilde{s}^2\left\{  {\tilde{s}\over 12}
-\langle {\tilde W}_-\varphi,\varphi \rangle\right\} >0.
\end{array}$$}
\end{proof}

The following Weitzenb\" ock formula will be used in the proof of the main result
which involves the Weyl curvature.

\beg{thm}[Weitzenb\" ock Formula \cite{bourguignon}]
On a Riemannian manifold, the Hodge/modern Laplacian can be
expressed in terms of the connection/rough Laplacian as
$$(d+d^*)^2=\nabla^*\nabla-2W+{s\over 3}$$
where $\nabla$ is the Riemannian connection and $W$ is the Weyl
curvature tensor. \end{thm}

\beg{thm}\label{finalthm} Let $M$ be a compact, 
Einstein-Hermitian 4-manifold of positive holomorphic bisectional curvature. Then the Betti number $b^2_-$ vanishes. \end{thm}

Here $b^2_-$ stands for the number of negative eigenvalues of the cup product in $H^2(M)$ as usual. By Hodge theory, this number is equal to the dimension of anti-self-dual harmonic $2$-forms on the Riemannian manifold $(M,g)$.

\beg{proof} These spaces are either K\"ahler, or else conformally K\" ahler by \cite{LeB:1995}. In the first case, according to the resolution \cite{SY:1980frankel} of Frankel's conjecture by Siu and Yau, the space has to be complex projective plane. By Theorem \ref{bgk}, the metric has to be the Fubini-Study metric. In the latter case, we can write $g=\tilde{s}^{-2}\tilde{g}$ for $g$ the Einstein, and $\tilde{g}$ the K\" ahler metric.

Assume, for a contradiction, there is a nonzero anti-self-dual $2$-form  $\varphi\in\Gamma(\Lambda^2_-)$ which is harmonic with respect to the K\"ahler metric $\tilde g$. Since we are on a compact manifold, we can rescale this form by a constant to have length strictly less than $1/\sqrt{2}$ everywhere. We write Weitzenb\" ock Formula for the K\" ahler metric $\tilde g$:
$$0=\Delta_{\tilde g} \varphi = \tilde\nabla^*\tilde\nabla \varphi -2\tilde W\varphi+{\tilde s\over 3}\varphi$$
Now we take the $L^2$-inner product (with respect to $\tilde g$) of both sides with $\varphi$:
\begin{align*}
0&= \int_M \langle\tilde\nabla\varphi,\tilde\nabla\varphi\rangle_{\tilde g} - 2\langle\tilde W\varphi,\varphi\rangle_{\tilde g} +{\tilde s\over 3}\langle\varphi,\varphi\rangle_{\tilde g} \, d\mu_{\tilde{g}}\\[10pt]
&=\int_M
|\tilde\nabla\varphi|^2_{\tilde{g}}+\tilde{s}^{-4}\left\{ {\tilde{s}\over 6}
-2\langle\tilde{W}_-\varphi,\varphi\rangle_g \right\} d\mu_{\tilde{g}}.
\end{align*}
At the points where $\varphi=0$, the term in the curly bracket is just the scalar curvature term, hence strictly positive. Let $p$ be a point where $\varphi|_p\neq 0$. Take an open set around $p$ on which $\varphi$ is nowhere zero.
Now working on this open set, let $\tilde{\varphi}:={\varphi/ |\varphi|\sqrt{2}}$ so that
$|\tilde{\varphi}|={1/\sqrt{2}}$. Since the conformally related $\tilde{\varphi}$ is also anti-self-dual and of constant norm,
the second estimate
${\tilde{s}\over 12}>\langle \tilde{W}_-\tilde{\varphi},\tilde{\varphi} \rangle$ is applicable for $\tilde{\varphi}$:
$$\large{
\begin{array}{rcl}
 \left\{ {\tilde{s}\over 6}
-2\langle\tilde{W}_-\varphi,\varphi\rangle_g \right\}_p
&=&{\tilde{s}\over 6}
-2\langle\tilde{W}_- \sqrt{2} |\varphi|\tilde{\varphi},\sqrt{2}|\varphi|\tilde{\varphi} \rangle_g\\[10pt]
&=&{\tilde{s}\over 6}
-4 |\varphi|^2 \langle\wt{W}_- \tilde{\varphi},\tilde{\varphi} \rangle_g\\[10pt]
&>&{\tilde{s}\over 6}
-4 .~{1\over2}~. ~{\tilde{s}\over 12}~=0 \end{array}}$$
at the point $p$.
Since the term in parenthesis is strictly positive everywhere we get $\nabla\varphi=0,\tilde{s}^{-4}=0$, a contradiction. So that no such $\varphi$ exists implying $b^2_-=0.$
\end{proof}

\begin{rmk}\textnormal{ Theorem \ref{finalthm}\label{florinremark} is actually true under the weaker assumption that \textit{orthogonal bisectional curvatures are positive}; that is,  $H(X,Y)>0$ whenever $X$ is perpendicular to $Y$ and $JY$. Indeed, in the non-K\"ahler case, when establishing the inequalities in the proof of Lemma \ref{firstestimate} and Lemma \ref{secondestimate}, we only considered the case when $\varphi = -\psi$, which corresponds to taking two \emph{orthogonal} complex lines $\mathrm{span}\{X,JX\}$ and $\mathrm{span}\{Y,JY\}$. On the other hand, in the K\"ahler-Einstein case, we notice that Theorem \ref{bgk} remains true if we weaken the assumption to positivity of orthogonal bisectional curvatures. See the proof of Theorem 5 in \cite{goldbergkobayashi67}, p.232, where the authors prove that the Einstein constant is a positive constant multiple of any of the orthogonal bisectional curvatures $H(X,JX)$, $H(X,Y)$ or $H(X,JY)$. We thank F. Belgun for this remark. }
\end{rmk}

\vspace{.05in}



{\small
\beg{flushleft}
\textsc{Orta mh. Z\"ubeyde Han\i m cd. No 5-3 Merkez 74100 Bart\i n
, T\" urk\'{i}ye.}\\
\textit{E-mail address:} \texttt{\textbf{kalafat@\,math.msu.edu}}
\end{flushleft}
}

{\small
\beg{flushleft}
\textsc{Department of Mathematics, NYC College of Technology, City University of New York, 300 Jay St, Brooklyn, NY 11201, USA.}\\
\textit{E-mail address:} \texttt{\textbf{ckoca@citytech.cuny.edu}}

\end{flushleft}
}


\bibliography{ehpbisec}{}

\begin{thebibliography}{EGH80}

\bibitem[And11]{anderssonthesis}
David Andersson.
\newblock {\em The 3-sphere and the bicycling black rings event horizon}.
\newblock Thesis. Stockholm University, 2011.

\bibitem[Ber65]{berger65}
Marcel Berger.
\newblock Sur quleques vari\'et\'es riemanniennes compactes d'{E}instein.
\newblock {\em C. R. Acad. Sci. Paris}, 260:1554--1557, 1965.

\bibitem[Bou81]{bourguignon}
Jean-Pierre Bourguignon.
\newblock Les vari\'et\'es de dimension {$4$} \`a signature non nulle dont la
  courbure est harmonique sont d'{E}instein.
\newblock {\em Invent. Math.}, 63(2):263--286, 1981.

\bibitem[Cal85]{C2:1985}
Eugenio Calabi.
\newblock Extremal {K}\"ahler metrics. {II}.
\newblock In {\em Differential geometry and complex analysis}, pages 95--114.
  Springer, Berlin, 1985.

\bibitem[Der83]{Derdzinski}
Andrzej Derdzi{\'n}ski.
\newblock Self-dual {K}\"ahler manifolds and {E}instein manifolds of dimension
  four.
\newblock {\em Compositio Math.}, 49(3):405--433, 1983.

\bibitem[DG91]{dragomirgrimaldi91}
Sorin Dragomir and Renata Grimaldi.
\newblock Cauchy-{R}iemann submanifolds of locally conformal {K}aehler
  manifolds. {III}.
\newblock {\em Serdica}, 17(1):3--14, 1991.

\bibitem[DO98]{dragomirornea}
Sorin Dragomir and Liviu Ornea.
\newblock {\em Locally conformal {K}\"ahler geometry}, volume 155 of {\em
  Progress in Mathematics}.
\newblock Birkh\"auser Inc., Boston, MA, 1998.

\bibitem[EGH80]{egh}
Tohru Eguchi, Peter~B. Gilkey, and Andrew~J. Hanson.
\newblock Gravitation, gauge theories and differential geometry.
\newblock {\em Phys. Rep.}, 66(6):213--393, 1980.

\bibitem[FIP04]{fip}
Maria Falcitelli, Stere Ianus, and Anna~Maria Pastore.
\newblock {\em Riemannian submersions and related topics}.
\newblock World Scientific Publishing Co., Inc., River Edge, NJ, 2004.

\bibitem[Fra61]{Fra:1961}
Theodore Frankel.
\newblock Manifolds with positive curvature.
\newblock {\em Pacific J. Math.}, 11:165--174, 1961.

\bibitem[GK67]{goldbergkobayashi67}
Samuel~I. Goldberg and Shoshichi Kobayashi.
\newblock Holomorphic bisectional curvature.
\newblock {\em J. Differential Geometry}, 1:225--233, 1967.

\bibitem[Hat02]{Hat:2002}
Allen Hatcher.
\newblock {\em Algebraic topology}.
\newblock Cambridge University Press, Cambridge, 2002.

\bibitem[KK15]{confk}
Mustafa Kalafat and Caner Koca.
\newblock Conformally {K}\"ahler surfaces and orthogonal holomorphic
  bisectional curvature.
\newblock {\em Geom. Dedicata}, 174:401--408, 2015.


\bibitem[KS16]{kalafatsari}
Mustafa Kalafat and Ramazan Sar\i.
\newblock On special submanifolds of the {P}age space.
\newblock {\em E-print available at
  \href{https://arxiv.org/abs/1608.03252}{arXiv:1608.03252.} math.DG}, 2016.
  
  
\bibitem[Koc12]{kocathesis}
Caner Koca.
\newblock {\em On {C}onformal {G}eometry of {K}\" ahler {S}urfaces}.
\newblock ProQuest LLC, Ann Arbor, MI. Ph.D. Thesis. Stony Brook University,
  2012.

\bibitem[Koc14]{kocapage}
Caner Koca.
\newblock Einstein {H}ermitian metrics of positive sectional curvature.
\newblock {\em Proc. Amer. Math. Soc.}, 142(6):2119--2122, 2014.


\bibitem[LeB97]{LeB:1995}
Claude LeBrun.
\newblock Einstein metrics on complex surfaces.
\newblock In {\em Geometry and physics ({A}arhus, 1995)}, volume 184 of {\em
  Lecture Notes in Pure and Appl. Math.}, pages 167--176. Dekker, New York,
  1997.

\bibitem[LeB12]{LeB:2011}
Claude LeBrun.
\newblock On {E}instein, {H}ermitian 4-manifolds.
\newblock {\em J. Differential Geom.}, 90(2):277--302, 2012.

\bibitem[Pag78]{P:1978}
Don~N. Page.
\newblock A {C}ompact {R}otating {G}ravitational {I}nstanton.
\newblock {\em Phys. Let.}, 79B(3):235--238, 1978.

\bibitem[Pag09]{page2arxiv}
Don~N. Page.
\newblock Some {G}ravitational {I}nstantons.
\newblock {\em E-print available at
  \href{https://arxiv.org/abs/0912.4922}{arXiv:0912.4922.} physics.gr-qc},
  2009.

\bibitem[SM99]{geofrust}
Jean-Fran\c{c}ois Sadoc and R\'emy Mosseri.
\newblock {\em Geometrical frustration}.
\newblock Collection Al\'ea-Saclay: Monographs and Texts in Statistical
  Physics. Cambridge University Press, Cambridge, 1999.

\bibitem[SY80]{SY:1980frankel}
Yum~Tong Siu and Shing~Tung Yau.
\newblock Compact {K}\"ahler manifolds of positive bisectional curvature.
\newblock {\em Invent. Math.}, 59(2):189--204, 1980.

\bibitem[Via11]{jeffdg}
Jeff Viaclovsky.
\newblock Lecture notes: Topics in riemannian geometry, 2011.
\newblock URL: \url{http://www.math.wisc.edu/~jeffv/courses/865_Fall_2011.pdf}.

\bibitem[Zha04]{zhangeulerforms}
Xiao Zhang.
\newblock Scalar flat metrics of {E}guchi-{H}anson type.
\newblock {\em Commun. Theor. Phys. (Beijing)}, 42(2):235--238, 2004.

\end{thebibliography}
\bibliographystyle{alphaurl}
\end{document}